\documentclass{elsarticle}
\usepackage{amssymb,amsbsy,amsmath,amsfonts}
\usepackage{graphicx} 
\usepackage{subfig}
\graphicspath{{figure/}}

\newcommand{\ZZ}{\mathbb{Z}}
\newcommand{\RR}{\mathbb{R}}
\newcommand{\NN}{\mathbb{N}}
\newcommand{\CC}{\mathbb{C}}

\newtheorem{theorem}{Theorem}[section]
\newtheorem{lem}[theorem]{Lemma}
\newtheorem{prop}[theorem]{Proposition}
\newtheorem{defn}[theorem]{Definition}
\newtheorem{rem}[theorem]{Remark}
\newtheorem{exm}[theorem]{Example}
\newtheorem{cor}[theorem]{Corollary}

\newenvironment{pf}
{\par\noindent{\bf Proof:} \hspace{0.03cm}}{{$\Box$} \par\medskip\par}

\journal{arXiv submission}

\begin{document}

\begin{frontmatter}

\title{Filters for anisotropic wavelet
  decompositions
}

\author[mc]{Mariantonia Cotronei}
\ead{mariantonia.cotronei@unirc.it}
\author[mr]{Milvia Rossini}
\ead{milvia.rossini@unimib.it}
\author[ts]{Tomas Sauer\corref{cor}}
\ead{Tomas.Sauer@uni-passau.de}
\author[mr]{Elena Volont\`e}
\ead{e.volonte2@campus.unimib.it}
\address[mc]{DIIES, Universit\`a Mediterranea di Reggio Calabria, Via
  Graziella - Feo di Vito, I-89122 Reggio Calabria, Italy} 
\address[ts]{Lehrstuhl f\"ur Mathematik mit Schwerpunkt Digitale
  Bildverarbeitung \& FORWISS, Universit\"at Passau, Fraunhofer IIS
  Research Group for Knowledge Bases Image Processing, Innstr. 43,
  D-94032 Passau, Germany}
\address[mr]{Universit\`a degli Studi di Milano--Bicocca, Via Cozzi
  55, I-20125 Milano, Italy}
\cortext[cor]{Corresponding author}


\begin{abstract}
  Like the continous shearlet transform and their relatives, discrete
  transformations based on the interplay between several filterbanks
  with \emph{anisotropic} dilations provide a high potential to
  recover directed features in two and more dimensions. Due to
  simplicity, most of the directional systems constructed so far were using
  prediction--correction methods based on interpolatory subdivision
  schemes. In this paper, we give a simple but effective construction
  for QMF (quadrature mirror filter) filterbanks which are the
  discrete object between orthogonal wavelet analysis. We also
  characterize when the filterbank gives rise to the existence of
  refinable functions and hence wavelets and give a generalized 
  shearlet construction for arbitrary dimensions and arbitrary
  scalings for which the filterbank construction ensures the existence
  of an orthogonal wavelet analysis.
\end{abstract}


\end{frontmatter}

\section{Introduction}
The problem of constructing \emph{multivariate} orthonormal wavelet
bases has been under investigation  for a long time, as many
applications, for example image or volume processing, require an
$s$-dimensional setting with $s>1$. 
The easiest way of proceeding is to consider ``separable'' bases,
which are simply tensor 
products of univariate wavelets. This corresponds to the case where
the underlying scaling matrix is diagonal. However, such an approach
has the drawback of  providing little flexibility in manipulating the
data since it privileges features aligned with the coordinate axis
directions. On the other hand, more recent approaches like shearlets
have been successful in even detecting \emph{directional}
singularities by means of anisotropic scaling.

The most general approach for such discrete methods is to consider any
expansive scaling matrix in $\ZZ^{s \times s}$, thus mapping the
integer lattice into an integer sub-lattice. However, in such a case, the
standard techniques used in the univariate setting cannot be used any more.
In fact, not only a Daubechies-like approach for constructing scaling
functions fails in more than one dimension, but  even in the case when
such an orthonormal scaling function is given, the construction of the
corresponding wavelets is far from straightforward and relies on
intricate algebraic principles, cf. \cite{Park95,park95:_algor_proof_of_suslin_stabil}.
 
In this paper we propose an approach to the construction of all the
filters associated with an orthonormal wavelet system which extends
tensor product constructions to arbitrary scaling matrices. It relies
on a Smith factorization of the scaling matrix, already used in
\cite{CotroneiSauerAl15} for the realization of bivariate
interpolatory  wavelet filters, and offers the possibility for any
easy extension of univariate filters associated to a general scaling
factor.
While this approach allows us to realize quadrature mirror filters
(QMF) for any expansive scaling matrix, the existence of the
corresponding scaling and wavelet \emph{functions}, defining a
multiresolution analysis, is more subtle and only holds with additional
assumptions on the scaling matrix. We will, however, prove a condition 
through the convergence of the associated subdivision scheme with
respect to a different scaling matrix.
 
The main motivation for such a construction is the realization of
multiple multiresolution analyses (MMRA), where several scaling
matrices are involved that possess anisotropic and  directional
components, thus making them useful in contexts where it is needed to
process data with directed lower dimensional features. The idea behind
an MMRA is that, at  each step of the filterbank decomposition,
different scaling matrices and filters are chosen from a finite
dictionary.
This approach generalizes the discrete shearlet transform from
\cite{KutyniokSauer09}, which considered only products of  parabolic
scaling and shears as scaling matrices. We specifically propose a
concept of a \emph{generalized shearlet system} with arbitrary
dilations for which it is easy to construct related orthogonal filter
banks that lead to an MMRA, making use of anisotropic diagonal
matrices and shears  of codimension $1$. 
We prove that such a choice provides \emph{slope resolution}, that is
the property of recovering hyperplanes of codimension $1$ by
applying appropriate combinations of the dilation matrices to a fixed
reference hyperplane.

\section{Notation and basic facts}
\label{sec:Notation}
We denote by $\ell (\ZZ^s)$ the space of all sequences, i.e., all functions from $\ZZ^s$
to $\RR$, and by $\ell_p (\ZZ^s)$ those sequences with finite
$p$--norm
$$
\| c \|_p = \left( \sum_{\alpha \in \ZZ^s} | c(\alpha) |^p
\right)^{1/p}, \qquad 1 \le p \le \infty.
$$
Extending these norms consistently to $p = 0,\infty$ as
$$
\| c \|_0 = \# \left\{ \alpha \in \ZZ^s : c(\alpha) \neq 0 \right\}
\qquad \text{and} \qquad
\| c \|_\infty = \sup_{\alpha \in \ZZ^s} | c(\alpha) |,
$$
we particularly use $\ell_{0} (\ZZ^s)$ for the subspace of  finitely
supported sequences.
To $a \in \ell_0 (\ZZ^s)$ we associate the \emph{symbol}
$$
a^\sharp (z) := \sum_{\alpha \in \ZZ^s} a(\alpha) \, z^\alpha, \qquad
z \in \CC_*^s := \left( \CC \setminus \{ 0 \} \right)^s
$$
which is a multivariate Laurent polynomial.
For a given integer matrix $\Theta \in \ZZ^{s\times s}$, we define
the \emph{dilation operator}  $D_{\Theta} : \ell(\ZZ^s) \to
\ell(\ZZ^s)$ by $c \mapsto c(\Theta \cdot)$.

An integer matrix $\Xi\in \ZZ^{s\times s}$ is called \emph{expansive}
or \emph{scaling matrix} if all of its eigenvalues are greater than
one in modulus, or, equivalently, of $\| \Xi^{-n} \| \to 0$ as $n \to
\infty$ in some matrix norm. Such a scaling matrix maps the integer
lattice $\ZZ^s$ to the sub-lattice $\Xi \ZZ^s$. It is well--known that
$\ZZ^s$ can be reconstructed as shifted sub-lattices,
$$
\ZZ^s = \bigcup_{\xi \in \ZZ^s_\Xi} \xi + \Xi \ZZ^s, \qquad \ZZ^s_\Xi =
\ZZ^s / \Xi \ZZ^s = \Xi [0,1)^s \cap \ZZ^s.
$$
The elements $\xi$ of $\ZZ^s_\Xi$ are called \emph{coset representers}
of $\ZZ^s$ modulo $\Xi$.

A function  $\phi$, which is, for example, assumed to belong to
$L_2(\RR)$, is said to be  \emph{refinable}  with respect to the
\emph{mask} $a\in \ell_{0} (\ZZ^s)$ and the scaling matrix $\Xi \in
\ZZ^{s \times s}$ if
$$
\phi=\sum_{\alpha\in \ZZ^s} a(\alpha)\phi(\Xi\cdot-\alpha).
$$
The function is called \emph{orthonormal} if it has orthonormal
integer translates, that is, 
\begin{equation}
  \label{eq:OrthonormDef}
  \langle \phi, \phi(\cdot -\gamma)
  \rangle:=\int_{\RR^s}\phi(x)\phi(x-\gamma)dx=\delta(\gamma),\qquad
  \gamma \in \ZZ^s,  
\end{equation}
where $\delta : \ZZ^s \to \RR$, $\delta (\gamma) = \delta_{\gamma,0}$
is usually called the \emph{pulse function} in signal processing.

Provided that \eqref{eq:OrthonormDef} holds, also
$\phi(\Xi^j\cdot)$, $j \in \ZZ$, has orthonormal integer translates
and, if we denote by
$$
V_j := \text{\rm span}\, \left\{ \phi ( \Xi^j \cdot - \alpha ) :
  \alpha \in \ZZ^s \right\}, \qquad j\in \NN,
$$
the space spanned by these translates, then $\{V_j\}_{j\in \NN}$
represents an orthogonal  \emph{multiresolution analysis} for the
space $L_2(\RR^s)$, cf. \cite{mallat09:_wavel_tour_signal_proces}.

The orthogonal complement of $V_j$ in $V_{j+1}$, denoted by $W_j =
V_{j+1} \ominus V_j$, is spanned by the translates and dilates of $d-1$
\emph{wavelets} $\psi_1,\dots, \psi_{d-1} \in V_1 \ominus V_0$, where
$d=\left| \det \Xi\right|$.
To label the wavelet functions, it is convenient to use the index set
$\ZZ^+_d :=\ZZ_d\setminus \{0\}$, where, as usual $\ZZ_d = \ZZ/d\ZZ:=
\{0,\dots,d - 1\}$. Belonging to $V_1$, the wavelets are all linear
combinations of translates and dilates of
the scaling function by means of the following relation:
$$
\psi_k=\sum_{\alpha\in \ZZ^s} b_k(\alpha)\phi(\Xi\cdot -\alpha),
\qquad b_k \in \ell_0 (\ZZ^s), \quad
k\in \ZZ^+_d, 
$$
and satisfy the orthogonality relations
$$
\langle \phi, \psi_\ell \rangle=0,\quad \langle \psi_k,
\psi_\ell(\cdot-\gamma) \rangle=\delta(\gamma)\delta_{k,\ell},\qquad
k,\ell \in \ZZ^+_d,\, \gamma \in \ZZ^s.
$$
An important property of a wavelet analysis are \emph{vanishing
	moments} which  guarantee, by a simple Taylor argument, fast decay of
the \emph{wavelet coefficients}
\begin{equation}
\label{eq:WaveletCoeff}
a_{k,\alpha}^j (f) := \frac{1}{|\det \Xi|^{j/2}} \int_{\RR^s}
f(x) \, \psi_k \left( \Xi^j x - \alpha \right) \, dx, \qquad j
\in \NN, \, k \in \ZZ_d^+, \, \alpha \in \ZZ^s,
\end{equation}
for smooth functions $f$. Here decay means decay with respect to the
level parameter $j$. Recall that the wavelets $\psi_k$, $k\in
\ZZ^+_d$, possess $n\ge 0$ \emph{vanishing moments} if  
$$
\int_{\RR^s} x^\beta \psi_k(x)\, dx=0, \qquad |\beta|\le n,\, k\in
\ZZ^+_d.
$$
The problem of constructing a scaling function for an MRA and,
consequently, a wavelet system, can be approached from the point of
view of subdivision schemes. Indeed, the limit function generated by a
convergent subdivision scheme is refinable and thus a candidate for
the scaling function spanning the spaces of an MRA. Conversely, any
scaling function must be the limit function of a convergent
subdivision scheme. Orthogonality of $\phi$ imposes further specific
conditions on the refinement mask. Moreover, the vanishing moment
property of the corresponding wavelets is connected to the polynomial
reproduction property of the scheme. We thus recall some basic
definitions on subdivision.

By $S_{\Xi,a} : \ell (\ZZ^s) \to \ell (\ZZ^s)$ we denote the
\emph{subdivision operator} with \emph{dilation matrix} $\Xi \in
\ZZ^{s \times s}$ and \emph{mask} $a \in \ell_{0} (\ZZ^s)$, defined as 
$$
S_{\Xi,a} c := \sum_{\alpha \in \ZZ^s} a( \cdot - \Xi \alpha ) \,
c(\alpha),
\qquad c \in \ell (\ZZ^s).
$$
Note that the finiteness
of the mask implies that $S_{\Xi,a}$ is a bounded operator from
$\ell_p (\ZZ^s)$ to itself for any $1 \le p \le \infty$.

\begin{defn}
  A \emph{subdivision scheme} is the iterative application of the
  operator $S_{\Xi,a}$, yielding, for any initial data $c^{[0]} \in \ell
  (\ZZ)$, a sequence
  $$
  c^{[r]} = S_{\Xi,a} \, c^{[r-1]} = S_{\Xi,a}^r \, c^{[0]}, \qquad r
  \in \NN.
  $$
  The subdivision scheme is called \emph{uniformly convergent} if
  for any $c \in \ell_\infty (\ZZ^s)$ there exists a uniformly
  continuous function $f_c \in C_u (\RR^s)$ such that
  $$
  \lim_{r \to \infty} \sup_{\alpha \in \ZZ^s} \left| f_c \left( \Xi^{-r}
      \alpha \right) -  S_{\Xi,a}^r \, c (\alpha) \right| = 0,
  $$
  and if the scheme is \emph{nontrivial}, i.e., there exists at least
  one $c \in \ell_\infty (\ZZ^s)$ such that $f_c \neq 0$.
\end{defn}

\begin{rem}
  The definition of convergence given here only considers the case $p
  = \infty$. There is, of course, convergence theory for $1 \le p <
  \infty$, see, for example
  \cite{DahmenMicchelli97,Jia95,MicchelliSauer97a}, but it is slightly
  more intricate though the core arguments are almost the same.
  Moreover, most of the popular wavelets are at least continuous. For
  the sake of simplicity, we thus restrict ourselves to uniform
  convergence here, extending the results to the $L_p$ case is straightforward.
\end{rem}

\begin{defn}
  We say that a subdivision operator $S_{\Xi,a}$ provides
  \emph{reproduction of degree $n$} or of \emph{order $n+1$} if
  \begin{equation}
    \label{eq:PolyPreserve}
    S_{\Xi,a} \Pi_k \subseteq \Pi_k, \qquad k=0,\dots,n,
  \end{equation}
  where
  $$
  \Pi_n := \mbox{\rm span}\, \left\{ (\cdot)^\alpha : |\alpha| \le n
  \right\},
  $$
  and the sequence $p \in \ell (\ZZ^s)$ attributed to a polynomial $p
  \in \Pi_n$ is simply $p |_{\ZZ^s}$.
\end{defn}

\begin{rem}
  Preservation of polynomials by subdivision operators is often
  described by different and contradicting terminology. It can mean,
  preservation of \emph{all} polynomials, i.e., $\left. S_{\Xi,a}
  \right|_{\Pi_n}$ is the identity, it can mean that the operator maps
  all $\Pi_k$, $k \le n$, to themselves, or that it only maps $\Pi_n$
  to itself. Our definition \eqref{eq:PolyPreserve} obviously means
  the second case.
\end{rem}

\noindent
Subdivision operators play a fundamental role in filterbanks,
cf. \cite{VetterliKovacevic95}. Recall
that a filterbank with dilation $\Xi$ and filters $a_j \in \ell_0
(\ZZ^s)$, $j \in \ZZ_n$, decomposes a signal $c \in \ell (\ZZ^s)$ into
its signal components 
$$
c_j := \downarrow_\Xi ( a_j * c ) := \left( a_j * c \right) (\Xi
\cdot), \qquad j \in \ZZ_n,
$$
by first filtering with $n$ (different) filters and then
\emph{downsampling} the signal. The mapping $F : c \to \left[ c_j : j \in
  \ZZ_n \right]$ is called the \emph{analysis filterbank}. To reverse
the process, the input signals $c_j$ are \emph{upsampled}, leading to
the nonzero entries
$$
\uparrow_\Xi c_j (\Xi \cdot) = c_j, \qquad j \in \ZZ_n.
$$
The upsampled data is then filtered by $b_j \in \ell_0 (\ZZ^s)$ and
the results are summed, i.e.,
$$
G : \left[ c_j : j \in \ZZ_n \right] \to \sum_{j \in \ZZ_n} b_j *
\left( \uparrow_\Xi c_j \right) = \sum_{j \in \ZZ_n} \sum_{\alpha \in
  \ZZ^s} b_j ( \cdot - \Xi \alpha ) \, c_j (\alpha) = \sum_{j \in \ZZ_n} S_{\Xi,b_j} c_j.
$$
In other words: \emph{reconstruction is subdivision}. Two concepts are
particularly important in the filterbank context.

\begin{defn}
  A filterbank is said to provide \emph{perfect reconstruction} if $FG
  = I$, and it is called \emph{critically sampled} if $n = |\det
  \Xi|$.
\end{defn}

\noindent
Most filterbanks in applications, especially those in the wavelet context, are
critically sampled perfect reconstruction filterbanks.

\section{Smith factorizations}
We recall that a matrix $\Theta \in \ZZ^{s \times s}$ is
called \emph{unimodular} if $| \det \Theta | = 1$ and that unimodular matrices
are exactly those integer matrices which also have an inverse in $\ZZ^{s
  \times s}$.

\begin{defn}
  Given $\Xi \in \ZZ^{s \times s}$, a decomposition of the form
  \begin{equation}
    \label{eq:XiSmith}
    \Xi = \Theta_1 \Sigma \Theta_2
  \end{equation}
  where $\Sigma$ is a diagonal matrix and $\Theta_1, \Theta_2$ are
  \emph{unimodular} matrices, is called a \emph{Smith factorization}
  of $\Xi$. \\
  The diagonal elements $\sigma_j$, $j=1,\dots,s$, of $\Sigma$ in the
  decomposition \eqref{eq:XiSmith} are called the \emph{Smith values}
  of this decomposition.
\end{defn}

\noindent
Smith factorizations can be computed, for example, by performing some
Gauss elimination with division by remainder and total pivoting to
diagonalize the matrix, cf. \cite{MarcusMinc69}. 
Note, however, that
Smith factorizations are not unique, neither with respect to the
diagonal matrix $\Sigma$, nor with respect to the unimodular factors
$\Theta_1$ and $\Theta_2$. In contrast to that, the \emph{Smith
  normal form} at least provides a unique $\Sigma$ by choosing
$\sigma_j$ as the quotient of the $j$-th and the
$j+1$-st \emph{determinantal divisor} of $\Xi$; recall that the $j$th
determinantal divisor is the greatest common divisor of all minors of
order $j$, i.e, the greatest common divisor (gcd) of all determinants
of $j \times j$ submatrices of $\Xi$. But even there the unimodular factors
$\Theta_1$ and $\Theta_2$ need not be unique at all.

\begin{exm}
  The diagonal matrix 
  $\Xi =
  \begin{bmatrix}
    3 \\ & 2 \\ & & 2 
  \end{bmatrix}$ has the nonzero minors $3$, $2$, $2$ of order $1$,
  the nonzero minors
  $$
  \det
  \begin{bmatrix}
    3 \\ & 2 
  \end{bmatrix} = 6, \qquad
  \begin{bmatrix}
    2 \\ & 2 
  \end{bmatrix} = 4
  $$
  of order $2$ and $\det \Xi = 12$ of order $3$. The greatest common
  divisors are $1$, $2$, and $12$, respectively and the values of the
  Smith normal form $1$, $2/1 = 2$ and $12/2 = 6$, hence the Smith
  normal form is
  $$
  \Xi = \Theta_1
  \begin{bmatrix}
    1 \\ & 2 \\ & & 6
  \end{bmatrix}
  \Theta_2.
  $$
\end{exm}

\noindent
Routines to compute the Smith normal form of a matrix are available,
for example in \texttt{matlab}.
There is an obvious relationship between matrices that have the same
Smith normal form which can be used to compute arbitrary Smith
factorizations. Let the normal form of $\Xi$ be given as
$\Xi = \Theta_1 \Sigma \Theta_2$ and let $\Sigma'$ be any diagonal
matrix with normal form $\Sigma' = \Lambda_1 \Sigma \Lambda_2$, then
$$
\Xi = \Theta_1 \Lambda_1^{-1} \Sigma' \Lambda_2^{-1} \Theta_2.
$$
This observation can be summarized as follows.

\begin{lem}\label{L:SmithEquiv}
  If $\Sigma$ is any diagonal matrix with the same Smith normal form
  as $\Xi$, then there exists a Smith decomposition
  $$
  \Xi = \Theta_1 \Sigma \Theta_2
  $$
  with proper unimodular matrices.
\end{lem}

\noindent
This lemma plays a role for example for scaling matrices like $\Sigma =
\begin{bmatrix}
  3 & \\ & 2 \\
\end{bmatrix}$ as considered in \cite{ElenaEtAlSauer2018S}. A Smith
normal form would always yield the diagonal $\begin{bmatrix}
  1 & \\ & 6 \\
\end{bmatrix}$, but by means of Lemma~\ref{L:SmithEquiv} there also
exist Smith factorizations with the diagonal $\Sigma$. Our later
construction will depend on choosing univariate refinable functions
whose scaling factors are the diagonals of $\Sigma$

\section{Orthogonal and quadrature mirror filters}\label{sec:qmf}
The following well--known fact is repeated for the reader's convenience.

\begin{lem}
  Let $\phi$ be a refinable function with respect to a mask $a$ and to the scaling matrix $\Xi$. If $\phi$ is orthonormal  then
  \begin{equation}
 \label{eq:aQMF}
 \sum_{\alpha \in \ZZ^s} a(\alpha) \, a (\alpha-\Xi
 \gamma ) = |\det \Xi| \, \delta ( \gamma ),  \quad \gamma \in \ZZ^s.
 \end{equation}
\end{lem}

\begin{pf}
   Straightforward computation yields for $\gamma \in \ZZ^s$
  \begin{eqnarray*}
  	\delta (\gamma)
  	& = & \int_{\RR^s} \phi (x) \, \phi( x-\gamma ) \, dx \\
  	& = & \sum_{\alpha,\beta \in \ZZ^s} \int_{\RR^s} a(\alpha) \phi
  	(\Xi x - \alpha) \, a(\beta) \phi (\Xi x - \Xi \gamma -
  	\beta) \, dx \\
  	& = &  \sum_{\alpha,\beta \in \ZZ^s} a(\alpha)
  	\, a(\beta) \int_{\RR^s} \phi (\Xi x - \alpha) \, 
  	\phi (\Xi x - \Xi \gamma - \beta) \, dx \\
  	& = & \frac{1}{|\det \Xi|} \sum_{\alpha \in \ZZ^s} a(\alpha) \, a(\alpha-
  	\Xi \gamma ),
  \end{eqnarray*}
  from which the claim follows immediately.
\end{pf}

\begin{rem}
	Note that the condition \eqref{eq:aQMF} means that the  auto-correlation
	$(a\star a)(-\Xi\cdot)$ is the
	mask of an \emph{interpolatory} subdivision scheme.
\end{rem}

\noindent
Based on this observation, we make the following definition.

\begin{defn}
  A mask $a \in \ell_{0} (\ZZ^s)$ is called \emph{orthonormal} with
  respect to the dilation matrix $\Xi$ if it satisfies the so-called
  \emph{QMF equation} \eqref{eq:aQMF}.
\end{defn}

\noindent
Following similar arguments as in  \cite{CotroneiSauerAl15}, we  now show how to construct
orthogonal filters by means of a generalized tensor product. To that end,
let  $\Xi \in \ZZ^{s \times s}$ be a dilation matrix with  decomposition of the form
(\ref{eq:XiSmith}) and Smith values $\sigma_j$, $j=1,\dots,s$.
Moreover, let $h_j$ be an univariate mask that satisfies a QMF equation with
scaling factor $\sigma_j$. We remark that it gives rise to a refinable
orthonormal limit function provided that the underlying subdivision
scheme converges. We set
$$
h := h_1 \otimes \cdots \otimes h_s, \qquad \text{i.e.,} \qquad
h(\alpha) = \prod_{j=1}^s h_j (\alpha_j), \quad \alpha \in \ZZ^s,
$$
and define
\begin{equation}
\label{eq:aOrthoDef}
a = h \left( \Theta_1^{-1} \cdot \right).
\end{equation}

\begin{lem}\label{L:XiQMF}
  The mask $a$ defined in \eqref{eq:aOrthoDef} satisfies the QMF
  equation \eqref{eq:aQMF}.
\end{lem}

\begin{pf}
  Since $h$ is, by construction, the mask of an orthogonal tensor
  product scheme, it satisfies $(h\star h) (- \Sigma \cdot ) = |\det
  \Sigma|\,  \delta$, from which we get that 
  \begin{eqnarray*}
    (a\star a)(-\Xi \cdot)
    & = & \sum_{\alpha \in \ZZ^s} a(\alpha) \, a(\alpha - \Theta_1 \Sigma \Theta_2
          \cdot )
          = \sum_{\alpha \in \ZZ^s} h \left( \Theta_1^{-1} \alpha
          \right) \, h \left(\Theta_1^{-1} \alpha- \Sigma \Theta_2 \cdot  
          \right) \\
    & = & \sum_{\alpha \in \ZZ^s} h \left( \alpha
          \right) \, h \left( \alpha -\Sigma \Theta_2 \cdot
          \right)
          = |\det \Sigma|\, \delta ( \Theta_2 \cdot ) = |\det
          \Xi|\,\delta (\cdot) 
  \end{eqnarray*}
  since $\Theta_2$ is unimodular, hence $\Theta_2 \gamma = 0$ iff
  $\gamma = 0$. This proves the claim.
\end{pf}

\noindent
In an analogous way, the orthogonal  wavelet filters associated to the
mask $a$ can be derived: let $g^k_j$, $k \in \ZZ_{\sigma_j}^+$, be the univariate
wavelet filters associated to $h_j$, $j=1,\dots,s$, 
and set, for consistency of notation, $g^{0}_j := h_j$. Then the
orthogonal wavelet filter system with
respect to the dilation factor $\Sigma$ is obtained via the following
tensor products 
\begin{equation}\label{eq:gtensor}
g_\eta := g^{\eta_1}_{1}\otimes g^{\eta_2}_{2}\otimes \dots \otimes
g^{\eta_s}_{s}, \qquad \eta \in \ZZ_\Sigma := \ZZ_{\sigma_1} \otimes
\cdots \otimes \ZZ_{\sigma_s},
\end{equation}
with $g_{(0,\ldots,0)}=h_1\otimes\dots \otimes h_s$. Like before, we finally set
\begin{equation}
  \label{eq:bOrthoDef}
  b_\eta := g_\eta \left( \Theta_1^{-1} \cdot \right), \qquad \eta \in \ZZ_\Sigma, 
\end{equation}
with  $b_{(0,\ldots,0)}=a$.
For convenience, we also write $\ZZ_\Sigma^+ := \ZZ_\Sigma \setminus \{
(0,\dots,0) \}$ to denote the index set for the high pass filters.

This construction automatically ensures that the complete system satisfies a filter
bank QMF identity. Indeed, as in the proof of Lemma~\ref{L:XiQMF} we
get for $\eta, \eta' \in \ZZ_\Sigma$, by taking into account $|\det
\Xi| = |\det \Sigma| = \prod_j |\sigma_j|$, that
\begin{eqnarray*}
  \lefteqn{
  \sum_{\alpha \in \ZZ^s} b_\eta (\alpha) \, b_{\eta'} \left( \alpha -
  \Xi \cdot \right) =
  \sum_{\alpha \in \ZZ^s} g_\eta \left( \Theta^{-1} \alpha \right) \,
  g_{\eta'} \left( \Theta^{-1} \alpha -
  \Sigma \Theta_2 \cdot \right)
  } \\
  & = & \sum_{\alpha \in \ZZ^s} g_\eta \left( \alpha \right) \,
        g_{\eta'} \left( \alpha - \Sigma \Theta_2 \cdot \right)
        = \prod_{j=1}^s \sum_{\alpha_j \in \ZZ} g_j^{\eta_j}
        (\alpha_j) \, g_j^{\eta_j'} \left( \alpha_j - \sigma_j \left(
        \Theta_2 \cdot \right)_j \right) \\
  & = & \prod_{j=1}^s \left| \sigma_j \right| \,
        \delta_{\eta_j,\eta_j'} \, \delta \left( \left( \Theta_2 \cdot
        \right)_j \right)
        = |\det \Xi| \, \delta_{\eta,\eta'} \, \delta (\cdot).
\end{eqnarray*}
This can be summarized in the following result.

\begin{theorem}\label{T:Filterbank}
  The filters $b_\eta$, $\eta \in \ZZ_\Sigma$, form a critically
  sampled QMF filterbank for the dilation matrix $\Xi$ with lowpass
  filter $b_0$ and highpass filters $b_\eta$, $\eta \neq 0$.
\end{theorem}

\begin{pf}
  We only have to verify the low and high pass properties of the
  filters which follow easily from
  $$
  b_\eta * 1 = \sum_{\alpha \in \ZZ^s} b_\eta (\alpha) = \prod_{j=1}^s
  \sum_{\alpha_j \in \ZZ} g_j^{\eta_j} (\alpha_j)
  $$
  and the respective high- and lowpass properties of the univariate filters.
\end{pf}

\noindent
If the mask $a$ defines a refinable function $\phi$, then the  \emph{wavelets} can be defined through the relation
\begin{equation}
\label{eq:WaveletDef}
\psi_\eta := \left( b_\eta * \phi \right) ( \Xi \cdot ) =
\sum_{\alpha \in \ZZ^s} b_\eta (\alpha) \, \phi (\Xi \cdot -
\alpha), \qquad \eta \in \ZZ_\Sigma.
\end{equation}
The orthogonality to $\phi$ is verified as usually via
\begin{eqnarray*}
	\lefteqn{\int_{\RR^s} \phi (x- \gamma) \, \psi_\eta (x) \, dx } \\
	& = & \sum_{\alpha,\beta \in \ZZ^s} a (\alpha) \, b_\eta (\beta)
	\int_{\RR^s} \phi (\Xi x - \Xi \gamma - \alpha) \, \phi (\Xi x -
	\beta) \, dx \\
	& = & \frac1{|\det \Xi|} \sum_{\alpha,\beta \in \ZZ^s} a (\alpha) \,
	b_\eta (\beta) \delta_{\alpha + \Xi \gamma,\beta}
	= \frac1{|\det \Xi|} \sum_{\alpha,\beta \in \ZZ^s} b_0 (\alpha
	- \Xi \gamma) \, b_\eta (\alpha) \\
	& = & \delta_{\eta,0} \, \delta (\gamma) = 0, \qquad \gamma \in
	\ZZ^s, \quad \eta \in \ZZ_\Sigma^+,
\end{eqnarray*}
and, by the same computation, via
$$
\int_{\RR^s} \psi_\eta (x- \gamma) \, \psi_{\eta'} (x) \, dx =
\delta_{\eta,\eta'} \delta(\gamma), \qquad \gamma \in \ZZ^s, \quad
\eta,\eta' \in \ZZ_\Sigma^+,
$$
we also get the orthonormality of the shifts of the wavelets
$\psi_\eta$.

By Theorem~\ref{T:Filterbank}, vanishing moments
are equivalent to preservation of polynomial spaces by the low pass
filter $b_0$, more precisely, by the associated \emph{subdivision
  operator} $S_{\Xi,b_0}$.

\begin{theorem}\label{T:VanishMom}
  If the univariate subdivision operators $S_{\sigma_j,g_j^0}$ provide
  reproduction of degree $n$, $j=1,\dots,s$, so does the operator
  $S_{\Xi,b_0}$ and thus the filterbank provides vanishing moments of
  degree $n$.
\end{theorem}

\begin{pf}
  By definition, we have that
  \begin{equation}
    \label{eq:VanishMomPf1}
    S_{\Xi,b_0} p = \sum_{\alpha \in \ZZ^s} b_0 ( \cdot - \Xi \alpha ) \,
    p (\alpha) = \sum_{\alpha \in \ZZ^s} g_0 ( \cdot - \Sigma \alpha ) \,
    p ( \Theta_2^{-1} \alpha).
  \end{equation}
  Since $p \in \Pi_n$ implies that $D_{\Theta_2^{-1}} p = p (
  \Theta_2^{-1} \alpha) \in \Pi_n$ as well, we can write it as
  $$
  p ( \Theta_2^{-1} \alpha) = \sum_{|\beta| \le n} p_{\beta} \,
  \alpha^\beta, \qquad \alpha \in \ZZ^s.
  $$
  Substituting this into \eqref{eq:VanishMomPf1}, we get that
  \begin{eqnarray*}
    S_{\Xi,b_0} p
    & = & \sum_{|\beta| \le n} p_{\beta} \sum_{\alpha \in \ZZ^s} g_0 (
          \cdot - \Sigma \alpha ) \, \alpha^\beta
          = \sum_{|\beta| \le n} p_{\beta} \sum_{\alpha \in \ZZ^s}
          \prod_{j=1}^s g^0_j \left( (\cdot)_j - \sigma_j \alpha_j
          \right) \alpha_j^{\beta_j} \\
    & = & \sum_{|\beta| \le n} p_{\beta} \, \prod_{j=1}^s \left(
          \sum_{\alpha_j \in \ZZ} g^0_j \left( (\cdot)_j - \sigma_j \alpha_j
          \right) \alpha_j^{\beta_j} \right)
          =: \sum_{|\beta| \le n} p_{\beta} \, \prod_{j=1}^s q_j (\alpha_j),
  \end{eqnarray*}
  with $q_j \in \Pi_{\beta_j}$, $j=1,\dots,s$, $|\beta| \le n$, by the
  assumption that the univariate schemes all reproduce $\Pi_n$. Set
  $q_\beta = \prod_j q_{\beta_j} \left( (\cdot)_j \right)$ and note
  that
  $$
  \deg q_\beta = \sum_{j=1}^s \deg q_{\beta_j} \le \sum_{j=1}^s
  \beta_j = |\beta|
  $$
  to finally conclude that $S_{\Xi,b_0} p$ is a polynomial and
  $$
  \deg S_{\Xi,b_0} p \le \deg \sum_{|\beta| \le n} p_{\beta} \,
  q_\beta
  \le \max_{|\beta| \le n} \deg q_\beta \le \max_{|\beta| \le n}
  |\beta| \le n.
  $$
  This verifies the polynomial reproduction property.
\end{pf}

\begin{rem}
  Observe that the tensor product scheme $S_{\Sigma,g_1^0 \otimes
    \cdots \otimes g_s^0}$ does not only preserve $\Pi_n$ under the
  assumption of Theorem~\ref{T:VanishMom}, but the much larger spaces
  $$
  \Pi_{\alpha} := \text{\rm span} \left\{ (\cdot)^\beta : \beta \le
    \alpha \right\}, \qquad \| \alpha \|_\infty \le n.
  $$
  This ``extra regularity'' gets lost, however, with the application
  of $\Theta^{-1}$.
\end{rem}

\section{Convergence and MRA generation}

We next want to see under which assumptions the filters $b_\eta$,
$\eta \in \ZZ_\Sigma$, are associated to a scaling function and
to corresponding wavelet functions, respectively, which generate an orthogonal
MRA. To this aim, we  study the convergence of the subdivision scheme
associated to the mask $a$.  

\noindent
While it is well--known that the tensor product subdivision scheme
$S_{\Sigma,h}$ converges whenever the coordinate schemes
$S_{\sigma_j,h_j}$, $j=1,\dots,s$, converge, nothing is known so far
about $S_{\Xi,a}$ where $a = h ( \Theta_1^{-1} \cdot)$. The answer is
given by the following observation.

\begin{prop}\label{P:EquivScheme}
  The subdivision scheme $S_{\Xi,a}$ converges if and only if
  $S_{\Sigma \Lambda,h}$ converges, with the unimodular matrix
  $\Lambda := \Theta_2 \Theta_1 \in \ZZ^{s \times s}$.
\end{prop}

\begin{pf}
  Like in \cite{CotroneiSauerAl15} we rewrite
  $$
  S_{\Xi,a} c = \sum_{\alpha \in \ZZ^s} h \left( \Theta_1^{-1} ( \cdot
    - \Theta_1 \Sigma \Theta_2 \alpha ) \right) c(\alpha)
  = \sum_{\alpha \in \ZZ^s} h \left( \Theta_1^{-1} \cdot
    - \Sigma \alpha \right) \, c( \Theta_2^{-1} \alpha)
  $$
  in terms of dilation operators
  as $S_{\Xi,a} = D_{\Theta_1^{-1}} S_{\Sigma,h} D_{\Theta_2^{-1}}$,
  which gives that
  $$
  S_{\Xi,a}^r = D_{\Theta_1^{-1}} \left( S_{\Sigma,h}
    D_{\Theta_2^{-1}} D_{\Theta_1^{-1}} \right)^r D_{\Theta_1}.
  $$
  Since
  $$
  D_{\Theta_2^{-1}} D_{\Theta_1^{-1}} c = \left( D_{\Theta_1^{-1}} c
  \right) \left( \Theta_2^{-1} \cdot \right)
  = c \left( \Theta_1^{-1} \Theta_2^{-1} \cdot \right) =
  D_{\Theta_1^{-1} \Theta_2^{-1}} c = D_{\Lambda^{-1}} (c),
  $$
  it follows from the unimodularity of $\Lambda$ that
  \begin{eqnarray*}
    S_{\Sigma,h} D_{\Theta_2^{-1}} D_{\Theta_1^{-1}} c
    & = & S_{\Sigma,h} D_{\Lambda^{-1}} c
          = \sum_{\alpha \in \ZZ^s} h \left(
          \cdot - \Sigma \alpha \right) \, c \left( \Lambda^{-1}
          \alpha \right) \\
    & = & \sum_{\alpha \in \ZZ^s} h \left(
          \cdot - \Sigma \Lambda \alpha \right) \, c(\alpha)
          = S_{\Sigma \Lambda,h} c,
  \end{eqnarray*}
  and therefore
  \begin{equation}
    \label{eq:LEquivSchemePf1}
    S_{\Xi,a}^r = D_{\Theta_1^{-1}} (S_{\Sigma \Lambda,h})^r D_{\Theta_1}.
  \end{equation}
  If we assume that $S_{\Xi,a}$ converges, then there exists a uniformly
  continuous function $\phi : \RR^s \to \RR$ such that
  $$
  \left( D_{\Theta_1} \left( S_{\Xi,a}^r \delta  - \phi ( \Xi^{-r}
    \cdot ) \right) \right) (\alpha)
  = S_{\Sigma \Lambda,h}^r \delta (\alpha) - \phi \left(
  \Xi^{-r} \Theta_1 \alpha \right)
  $$
  tends to zero as $r \to \infty$ and this convergence occurs
  uniformly in $\alpha$. Observing that
  $$
  \Xi^{-r} \Theta_1 = \left( \Theta_2^{-1} \Sigma^{-1} \Theta_1^{-1}
  \right)^r \Theta_1 = \Theta_1 \left( \Theta_1^{-1} \Theta_2^{-1}
    \Sigma \right)^r = \Theta_1 \left( \Sigma \Lambda \right)^{-r}
  $$
  then allows us to conclude, together with the unimodularity of
  $\Theta_1$, that
  $$
  \lim_{r \to \infty} \sup_{\alpha \in \ZZ^s} \left| S_{\Sigma
      \Lambda,h}^r \delta - \phi \left( (\Sigma \Lambda)^{-r} \cdot
    \right) \right| (\alpha) = \lim_{r \to \infty} \sup_{\alpha \in \ZZ^s}
  \left| S_{\Xi,a}^r \delta  - \phi ( \Xi^{-r}
    \cdot ) \right| (\alpha) = 0,
  $$
  hence
  $\tilde\phi := \phi \left( \Theta_1 \cdot
  \right)$ is a \emph{basic limit function} for $S_{\Sigma \Lambda,h}$. The
  converse is obtained in the same way by realizing that the limit
  $\eta$ for the right hand side defines a limit $\phi := \tilde \phi (
  \Theta_1^{-1} \cdot )$ for the left hand side.
\end{pf}

\begin{rem}
  Convergence of subdivision schemes can also be considered on $L_p$
  spaces, $1 \le p < \infty$,
  cf. \cite{DahmenMicchelli97,Jia95,MicchelliSauer98}, even for
  arbitrary dilation matrices. It is straightforward though tedious to
  adapt the above arguments to arbitrary $L_p$ convergence. 
\end{rem}

\noindent
Given a convergent scheme $S_{\Sigma,h}$, the convergence of
$S_{\Sigma \Lambda,h}$ is not easy to check. There is however a
situation in which it is obvious.

\begin{cor}\label{C:SimilarConv}
  If $\Xi$ is \emph{similar} to $\Sigma$ over $\ZZ^{s \times s}$, that
  is, $\Xi = \Theta \Sigma \Theta^{-1}$ for some $\Theta \in \ZZ^{s
    \times s}$, then the convergence of $S_{\Xi,a}$ and $S_{\Sigma
    \Lambda,h}$ are equivalent. 
\end{cor}

\begin{pf}
  Just note that $\Lambda = \Theta^{-1} \Theta = I$.
\end{pf}

By collecting all the previous results and remarks, we are now able to state the main results of this section.

\begin{theorem}\label{T:PsiAreWavelets}	
	Let  $\Xi \in \ZZ^{s \times s}$ be a dilation matrix with  decomposition of the form
	(\ref{eq:XiSmith}), with $\Theta_1^{-1} \Theta_2 = I$.
	 Let $\tilde \phi=\otimes_{j=1}^s\phi_{\sigma_j}$ be the tensor product of univariate orthogonal scaling functions associated to the dilations $\sigma_j$, $j=1,\dots,s$, representing the Smith values of $\Xi$. Then the function $\phi=\tilde \phi (\Theta_1^{-1}\cdot)$ generates an orthogonal MRA in $L_2(\RR^s)$  with respect to the dilation $\Xi$ and the functions defined as in 	(\ref{eq:WaveletDef})
	are the corresponding orthonormal wavelets.
	\end{theorem}
\begin{pf}
	The statement is a consequence of all the previous results. In particular, from  Theorem \ref{P:EquivScheme} and Corollary \ref{C:SimilarConv} it follows that $\phi$ generates an orthogonal  MRA and is refinable since it is associated to a convergent subdivision scheme. Furthermore  the functions defined in  (\ref{eq:WaveletDef}) generate the complementary spaces of such MRA and  satisfy the orthonormality conditions, as proved in Section \ref{sec:qmf}, thus are the associated wavelets.
\end{pf}

We conclude the section by underlining that the previous result
guarantees the construction of an MRA associated to a generic dilation
$\Xi$, provided that it is similar to a diagonal matrix, starting from
orthonormal univariate scaling functions of any scaling factor. The
existence and construction of such Daubechies-like functions have been
extensively studied, we just mention
\cite{chui95:_const,heller95:_rank_m_n} here. And
even in the case when for a general dilation convergence of the
subdivision scheme cannot be assured, it is still  possible, with our
proposed construction, to realize an orthogonal filterbank. In
particular, such a filterbank provides perfect reconstruction and
vanishing moments and is therefore suitable from a purely signal
processing perspective.

\section{Multiple multiresolution and a generalized shearlet system}
One main motivation for wavelet schemes with arbitrary dilation
matrices comes from \emph{multiple multiresolution analysis} (MMRA)
and the extension of the discrete shearlet transform from
\cite{KutyniokSauer09}. Recall from
\cite{Sauer11:_shear_multir_multip_refin} that an MMRA is based on $m
\ge 1$ filterbanks with dilations $\Xi_j = \Theta_{j1} \Sigma_j
\Theta_{j2}$, $j \in \ZZ_m$, and filters $B_j = \left( b^j_\eta : \eta
  \in \ZZ_{\Sigma_j} \right)$. The (minimal) assumptions to be made are that
\begin{enumerate}
\item the dilation matrices are \emph{jointly expanding} which is best
  described by $\rho \left( \Xi_j^{-1} : j \in \ZZ_m \right) < 1$,
  where $\rho$ denotes the \emph{joint spectral radius} of the
  dilations.
\item the subdivision scheme $S_{\Xi_0,b_0^0}$ corresponding to the
  low pass synthesis filter $b_0^0$ converges.
\end{enumerate}
Writing
$$
\ZZ_m^* := \bigcup_{n=1}^\infty \ZZ_m^n
$$
for the set of all \emph{finite} sequences in $\ZZ_m$, we can
associate to any $\mu \in \ZZ_m^*$  the function
$$
\phi_\mu := \lim_{r \to \infty} S_{\Xi_0,b_0^0}^r
S_{\Xi_{d_n},b_0^{d_n}} \cdots S_{\Xi_{d_1},b_0^{d_1}} \, \delta.
$$
These functions satisfy the \emph{joint refinement equation}
\begin{equation}
  \label{eq:JointRefinement}
  \phi_{(j,\mu)} = \sum_{\alpha \in \ZZ^s} b_0^j (\alpha) \, \phi_\mu
  \left( \Xi_j \cdot - \alpha \right), \qquad j \in \ZZ_m, \, \mu \in
  \ZZ_m^*, 
\end{equation}
which allows us to build a redundant multiresolution,
cf. \cite{Sauer11:_shear_multir_multip_refin} for details.

From a filterbank perspective an MMRA recursively applies the analysis
filterbanks $F_j$ based on $\Xi_j$ and $B_j$, $j \in \ZZ_m$, to the
initial data and iteratively applies these decompositions on the low
pass components, leading to a \emph{tree} of wavelet
decompositions. The intuition, originating from \cite{KutyniokSauer09}
is that each branch of the tree detects features with a certain
directional component that is encoded in the digit sequence defining
the branch. \emph{Slope resolution}, the property that we will verify
in Theorem~\ref{T:SlopeApproximation}, then ensures that in the limit
all possible directions are met with arbitrary accuracy by some branch
in the tree.

From an analysis point of view, the spaces of the MMRA are now defined
in the highly redundant form
\begin{equation}
  \label{eq:MMRASpace}
  V_j = \text{\rm span}\, \left\{ \phi_\mu \left( \Xi_\eta \cdot -
      \alpha \right) : \mu \in \ZZ_m^*, \, \eta \in \ZZ_m^j, \, \alpha
    \in \ZZ^s \right\}
\end{equation}
as dilates and shifts of \emph{all} possible limit functions in the
process. The orthogonality concept in the multiresolution, on the
other hand, works \emph{along the trees}. More precisely, we have that
$$
\left\langle \phi_\mu, \psi_{\eta,(\mu',\mu)} (\cdot - \alpha) \right\rangle = 0,
\qquad \mu,\mu' \in \ZZ_m^*, \, \eta \in \ZZ_{\Sigma_j}^+, \, \alpha \in \ZZ^s,
$$
where the proof is as in \cite{Sauer11:_shear_multir_multip_refin}.

Together with the results from the preceding section, it is now easy
to construct such a multiple multiresolution provided that we know a
wavelet scheme for $\Xi_0$. Constructing the remaining $b_\eta^j$
according to the preceding recipe then leads to an MMRA which inherits
orthogonality as well as vanishing moments from the univariate
schemes. However, this approach has one drawback: it is not symmetric
with respect to the dilations by giving $\Xi_0$ a very special
role. In the special case, however, that all dilation matrices are
unimodularly similar to $\Xi_0$, we can give a more symmetric
construction generalizing the discrete shearlets, which we will point
out in the remainder of this section.

Based on the observation of Corollary~\ref{C:SimilarConv}, we now
define a concept of generalized shearlets with arbitrary dilations for
which it is easy to construct related orthogonal filter banks that
lead to a convergent multiple subdivision scheme. In particular, it
will not be
necessary to use special properties of the parabolic shearing as in
\cite{KutyniokSauer09}.
The shears we are using are those of
codimension $1$. That is, with the canonical unit vectors $e_j \in
\RR^{s-1}$, we set
\begin{equation}
  \label{eq:GammajDef}
  \Gamma_j :=
  \begin{bmatrix}
    I_{s-1} & -e_j \\ & 1
  \end{bmatrix}, \qquad j=1,\dots,s-1,  
\end{equation}
and extend it to $\Gamma_0 = I$, i.e., $e_0 = 0$. For two nonnegative
integers $\sigma_1 > \sigma_2 > 1$ we then set
\begin{equation}
 \label{eq:Xij}
\Xi_0 :=
\begin{bmatrix}
  \sigma_1 \, I & \\ & \sigma_2
\end{bmatrix},
\qquad \Xi_j := \Gamma_j^{-1} \Xi_0 \Gamma_j =
\begin{bmatrix}
  \sigma_1 I & \left( \sigma_2 - \sigma_1 \right) e_j \\
  & \sigma_2
\end{bmatrix},
\end{equation}
which implies that
\begin{equation}
  \label{eq:XijInv}
  \Xi_j^{-1} =
  \begin{bmatrix}
    \sigma_1^{-1} I & \frac{\sigma_1 - \sigma_2}{\sigma_1 \, \sigma_2}
    \, e_j \\
    & \sigma_2^{-1}
  \end{bmatrix}.  
\end{equation}
In multiple subdivision,
cf. \cite{Sauer11:_shear_multir_multip_refin,Sauer11:_multip_subdiv},
which is the basis of shearlet multiresolution systems, one considers
the iterated matrices
\begin{equation}
  \label{eq:XepsDef}
  \Xi_\epsilon = \Xi_{\epsilon_n} \cdots \Xi_{\epsilon_1}, \qquad
  \epsilon \in \ZZ_s^n, \quad n \in \NN.
\end{equation}
We can give an explicit expression for these matrices and their
inverses.

\begin{lem}
  For $n \in \NN$ and $\epsilon \in \ZZ_s^n$ we have that
  \begin{equation}
    \label{eq:XiepsInvFormula}
    \Xi_\epsilon^{-1} =
    \begin{bmatrix}
      \sigma_1^{-n} I & \sigma_2^{-n} p_\epsilon \left(
        \frac{\sigma_2}{\sigma_1} \right) \\
      & \sigma_2^{-n}
    \end{bmatrix}, \qquad p_\epsilon (x) = (1-x) \sum_{k=1}^n x^{k-1} e_{\epsilon_k},
  \end{equation}
  and 
  \begin{equation}
    \label{eq:XiepsFormula}
    \Xi_\epsilon =
    \begin{bmatrix}
      \sigma_1^{n} I & -\sigma_1^{n} p_\epsilon \left(
        \frac{\sigma_2}{\sigma_1} \right) \\
      & \sigma_2^{n}
    \end{bmatrix}.
  \end{equation}
\end{lem}

\begin{pf}
  We use induction on $n$, where the case $n=1$ simply consists of
  rewriting the upper right entry of \eqref{eq:XijInv} as
  $$
  \sigma_2^{-1} \left( 1 - \frac{\sigma_2}{\sigma_1} \right)
  e_{\epsilon_1} = \sigma_2^{-1} \, p_\epsilon \left(
    \frac{\sigma_2}{\sigma_1} \right), \qquad \epsilon \in \ZZ_s^1.
  $$
  To advance the induction hypothesis we compute, for $\epsilon =
  (\hat \epsilon,\epsilon_{n+1}) \in \ZZ_s^{n+1}$,
  \begin{eqnarray*}
    \Xi_\epsilon^{-1}
    & = & \Xi_{\hat \epsilon}^{-1} \Xi_{\epsilon_{n+1}}
          = \begin{bmatrix}
      \sigma_1^{-n} I & \sigma_2^{-n} p_{\hat \epsilon} \left(
        \frac{\sigma_2}{\sigma_1} \right) \\
      & \sigma_2^{-n}
    \end{bmatrix} \begin{bmatrix}
      \sigma_1^{-1} I & \frac{\sigma_1 - \sigma_2}{\sigma_1 \, \sigma_2}
      \, e_{\epsilon_{n+1}} \\
      & \sigma_2^{-1}
    \end{bmatrix} \\
    & = & \begin{bmatrix}
      \sigma_1^{-n-1} I & \sigma_1^{-n} \frac{\sigma_1 -
        \sigma_2}{\sigma_1 \, \sigma_2} 
      \, e_{\epsilon_{n+1}}+ \sigma_2^{-n-1} p_{\hat \epsilon} \left(
        \frac{\sigma_2}{\sigma_1} \right) \\
      & \sigma_2^{-n-1}
    \end{bmatrix} 
  \end{eqnarray*}
  and the expression in the upper right simplifies after inserting the
  induction hypothesis as follows
  \begin{eqnarray*}
    \lefteqn{
    \sigma_1^{-n-1} \sigma_2^{-n-1} \left( \sigma_1
    \sigma_2^n - \sigma_2^{n+1} \right)
    e_{\epsilon_{n+1}} + \sigma_2^{-n-1} \left(
    1 - \frac{\sigma_2}{\sigma_1} \right) \sum_{k=1}^n \left(
    \frac{\sigma_2}{\sigma_1} \right)^{k-1} e_{\epsilon_k}
    }\\
    & = & \sigma_2^{-n-1} \left( \left( \left(
          \frac{\sigma_2}{\sigma_1} \right)^n - \left(
          \frac{\sigma_2}{\sigma_1} \right)^{n+1} \right)
          e_{\epsilon_{n+1}} + \left( 1 - \frac{\sigma_2}{\sigma_1} \right)
          \sum_{k=1}^n \left(\frac{\sigma_2}{\sigma_1} \right)^{k-1}
          e_{\epsilon_k} \right) \\
    & = & \sigma_2^{-n-1}  \left( 1 - \frac{\sigma_2}{\sigma_1} 
          \right) \left( \sum_{k=1}^n \left(\frac{\sigma_2}{\sigma_1} \right)^{k-1}
          e_{\epsilon_k} + \left(\frac{\sigma_2}{\sigma_1}
          \right)^n e_{\epsilon_{n+1}} \right) \\
    & = & \sigma_2^{-n-1} p_\epsilon \left( \frac{\sigma_2}{\sigma_1}
          \right).
  \end{eqnarray*}
  Rewriting this as
  $$
  \Xi_\epsilon^{-1} =
  \begin{bmatrix}
    I & p_\epsilon \left( \frac{\sigma_2}{\sigma_1} \right) \\
    & 1
  \end{bmatrix}
  \begin{bmatrix}
    \sigma_1^{-n} I & \\
    & \sigma_2^{-n}
  \end{bmatrix},
  $$
  equation \eqref{eq:XiepsFormula} follows readily.
\end{pf}

\noindent
Using \eqref{eq:XiepsInvFormula}, we find that
$$
\| \Xi_\epsilon^{-1} \|_1 \le \max \left\{ \sigma_2^{-n},
  \sigma_1^{-n} + \sigma_2^{-n} \, \left| p_\epsilon \left(
      \frac{\sigma_2}{\sigma_1} \right) \right| \right\},
$$
and since the second expression is bounded by
\begin{eqnarray*}
  \lefteqn{\sigma_2^{-n} \left( \left( \frac{\sigma_2}{\sigma_1}
  \right)^n + \left( 1 - \frac{\sigma_2}{\sigma_1} \right) \sum_{k=1}^n
  \left( \frac{\sigma_2}{\sigma_1} \right)^{k-1} \right)
  } \\
  & \le & \sigma_2^{-n} \left( \left( \frac{\sigma_2}{\sigma_1}
          \right)^n + \left( 1 - \frac{\sigma_2}{\sigma_1} \right)
          \sum_{k=1}^\infty \left( \frac{\sigma_2}{\sigma_1}
          \right)^{k-1} \right) 
          = \sigma_2^{-n} \left( 1 + \left( \frac{\sigma_2}{\sigma_1}
          \right)^n \right),
\end{eqnarray*}
it follows that
$$
\left\| \Xi_\epsilon^{-1} \right\|_1^{1/n} < \frac{\max \left\{ 1,\left(
      1+ \left( \frac{\sigma_2}{\sigma_1} \right)^n \right)^{1/n}
  \right\}}{\sigma_2}
$$
which is $< 1$ for any $\epsilon \in \ZZ_s^n$ provided that $n$ is
sufficiently large since $\left( 1+ \left( \frac{\sigma_2}{\sigma_1}
  \right)^n \right)^{1/n} \to 1$ and $\sigma_2 > 1$. In other words,
the matrices $\Xi_j$ are \emph{jointly contractive}.

The last property of the shearlet analysis we need to prove is
\emph{slope resolution}, cf. \cite{CotroneiSauerAl15} which is the
property to recover directions 
from basic directions by means of an appropriate $\Xi_\epsilon$. To
that end, we associate to any \emph{hyperplane}
$$
H = H_v := \left\{ x \in \RR^s : v^T x = 0 \right\}, \qquad v \in
\RR^s,
$$
the \emph{slope} $w \in \RR^{s-1}$ if $v = (w,1)^T$. Any normal
direction $v$ with $v_s = 0$ is excluded here, and we will comment on
that later. Finally, we define the $k$--dimensional standard simplex as
$$
\mathbb{S}_k := \left\{ x \in \RR^k : x_j \ge 0, \sum_{j=1}^k x_j \le
  1 \right\} \subset \RR^k.
$$
Now we have the following result.

\begin{theorem}\label{T:SlopeApproximation}
  For any $w \in \mathbb{S}_{s-1}$, any $w' \in \RR^{s-1}$ and any
  $\delta > 0$ there exist $n \in \NN$ and $\epsilon \in \ZZ_s^n$ such
  that
  \begin{equation}
    \label{eq:SlopeApproximation}
    \left\|
      \begin{bmatrix}
        w' \\ 1 \\
      \end{bmatrix}
      - \sigma_2^{n} \Xi_\epsilon^{-1} \,
      \begin{bmatrix}
        w \\ 1 \\
      \end{bmatrix}
    \right\| < \delta.
  \end{equation}
\end{theorem}

\noindent
Before we prove this result, let us briefly recall its geometric
meaning. The vectors $v = (w,1)^T$ and $v' = (w',1)^T$ can be seen as
normals of two hyperplanes $H$ and $H'$. The estimate
\eqref{eq:SlopeApproximation} now says that the hyperplane $H'$,
corresponding to a directed singularity at some point, can be obtained
by applying $\Xi_\epsilon$ to the \emph{reference hyperplane} $H$, so
that all possible directions, except those with last component equal
to zero, can be constructed in the associated multiple multiresolution
analysis, just like in the shearlet case.

\begin{pf}
  The proof is a slight generalization and modification of the one
  given in
  \cite{ElenaEtAlSauer2018S}. We first note that, by \eqref{eq:XijInv},
  $$
  \Xi_j^{-1}
  \begin{bmatrix}
    w \\ 1 \\
  \end{bmatrix}
  = 
  \begin{bmatrix}
    \sigma_1^{-1} I & \frac{\sigma_1 - \sigma_2}{\sigma_1 \, \sigma_2}
    \, e_j \\
    & \sigma_2^{-1}
  \end{bmatrix} \, \begin{bmatrix}
    w \\ 1 \\
  \end{bmatrix}
  = \frac{1}{\sigma_2}
  \begin{bmatrix}
    \frac{\sigma_2}{\sigma_1} w + \left( 1 - \frac{\sigma_2}{\sigma_1}
    \right) e_j \\ 1
  \end{bmatrix},
  $$
  which induces the contractions $h_j (w) := \frac{\sigma_2}{\sigma_1}
  w + \left( 1- \frac{\sigma_2}{\sigma_1} \right) e_j$,
  $j\in \ZZ_s$, on $\RR^{s-1}$, that satisfy
  $$
  h_j \left( \mathbb{S}_{s-1} \right) =
  \frac{\sigma_2}{\sigma_1} \mathbb{S}_{s-1} + \left(
    1 - \frac{\sigma_2}{\sigma_1} \right) \subset
  \mathbb{S}_{s-1}, \qquad j \in \ZZ_s,
  $$
  and
  $$
  \mathbb{S}_{s-1} = \bigcup_{j \in \ZZ_s} h_j \left( \mathbb{S}_{s-1}
  \right),
  $$
  so that $\mathbb{S}_{s-1}$ is an invariant set for the $h_j$ which
  yields that for any compact $X \subset \RR^{s-1}$ we have
  \begin{equation}
    \label{eq:TSlopeApproximationPf1}
    \mathbb{S}_{s-1} = \lim_{n \to \infty} \bigcup_{\epsilon \in
      \ZZ_s^n} h_\epsilon (X), \qquad h_\epsilon := h_{\epsilon_n}
    \circ \cdots \circ h_{\epsilon_1},
  \end{equation}
  in the Hausdorff norm, see \cite{Hutchinson81}.
  Since, by \eqref{eq:XiepsInvFormula},
  $$
  \Xi_\epsilon^{-1}
  \begin{bmatrix}
    w \\ 1 \\
  \end{bmatrix} = \Xi_{\epsilon_n}^{-1} \, \Xi_{\epsilon_{n-1}}^{-1}
  \cdots \Xi_{\epsilon_1}^{-1}  \begin{bmatrix}
    w \\ 1 \\
  \end{bmatrix},
  $$
  it follows by induction that
  $$
  \Xi_\epsilon^{-1}
  \begin{bmatrix}
    w \\ 1 \\
  \end{bmatrix} = \sigma_2^{-n} 
  \begin{bmatrix}
    \left( \frac{\sigma_2}{\sigma_1} \right)^n \, w + p_{\epsilon}
    \left( \frac{\sigma_2}{\sigma_1} \right) \\ 1 
  \end{bmatrix} = \sigma_2^{-n} 
  \begin{bmatrix}
    h_{\epsilon} (w) \\ 1 
  \end{bmatrix}.
  $$
  Now choose $w' \in -\mathbb{S}_{s-1}$ and $\delta > 0$, then there
  exists by \eqref{eq:TSlopeApproximationPf1} with $X = \{ w \}$, an
  index $\epsilon \in \ZZ_s^n$ for some $n \in \NN$, such that
  $$
  \| w' - h_\epsilon (w) \| < \delta.
  $$
  Thus,
  $$
  \delta > \left\|
    \begin{bmatrix}
      w' \\ 1 \\
    \end{bmatrix} - \begin{bmatrix}
      h_\epsilon (w) \\ 1 \\
    \end{bmatrix}
  \right\| = \left\|
    \begin{bmatrix}
      w' \\ 1 \\
    \end{bmatrix} - \sigma_2^n \, \Xi_{\epsilon}^{-1} \, \begin{bmatrix}
      w \\ 1 \\
    \end{bmatrix}
  \right\|
  $$
  as claimed.
\end{pf}

\noindent
We have to recall an important fact concerning the true complexity of
discrete shearlet methods.

\begin{rem}\label{R:SlopeApproximation}
  Theorem~\ref{T:SlopeApproximation} only yields an approximation
  result for slopes in $-\mathbb{S}_{s-1}$. To obtain slopes in
  $$
  \mathbb{S}_{s-1}^\eta = \left\{ x \in \RR^{s-1} : \eta_j \, x_j \ge
    0, \sum_{j \in \ZZ_s} \eta_j x_j \le 1 \right\}, \qquad \eta \in
  \{0,1\}^{s-1},
  $$
  we have to replace $\Gamma_j$ by
  $$
  \Gamma_j^\eta :=
  \begin{bmatrix}
    \sigma_1 \, I & (-1)^{\eta_j+1} \, e_j \\ & \sigma_2
  \end{bmatrix},
  $$
  yielding a different multiresolution for any sign distribution $\eta
  \in \{0,1\}^{s-1}$. This is a known drawback of multivariate discrete
  shearlet systems that a full directional resolutions requires
  $2^{s-1}$ of them to be run in parallel, namely, the ones based on
  the systems $\Gamma_j^\eta$, $\eta \in \{0,1\}^s$.

  And even that is not enough: to capture the hyperplanes with normals
  $e_j$ requires the same to be run with an anisotropic scaling matrix
  of the form $\sigma_2 I + (\sigma_1 - \sigma_2) e_j e_j^T$ that has
  $\sigma_2$ in the $j$-th component. In summary, one needs a total of
  $s \, 2^{s-1}$ multiple multiresolution analyses. Keep in mind that
  this fact has already been observed for shearlets in
  \cite{Sauer11:_shear_multir_multip_refin}.
\end{rem}

\noindent
Let us summarize the achievements of this section. For the arbitrary
dilations of codimension $1$ from \eqref{eq:Xij}, the construction of
section \ref{sec:qmf} gives
us orthogonal filters whose associated subdivision schemes \emph{all}
converge. Since they all satisfy
$$
\Sigma_j \Lambda_j = \Sigma_j = \Xi_0,
$$
it even follows that that the associated \emph{multiple subdivision scheme}
converges, cf. \cite{Sauer11:_multip_subdiv}. Moreover, the filters
give rise to an orthogonal multiple multiresolution analysis with
slope resolution for which scaling functions and wavelets exist. The
approach includes parabolic scalings where $\sigma_1 = \sigma_2^2$ and
$p_\epsilon$ is evaluated at integers which can be seen as
$\sigma_2$--adic digit expansion of the slope,
cf. \cite{KutyniokSauer09,Sauer11:_shear_multir_multip_refin}. On the
other hand, Theorem~\ref{T:SlopeApproximation} shows that the 
evaluation of $p_\epsilon$ at arbitrary \emph{rational points} in
$(1,\infty)$ also gives a unique relationship between sloped and digits.
This again suggests the multiresolution based on $\sigma_1 = 3$, and
$\sigma_2 = 2$ with the smallest expansive integer factors on the
diagonal as considered first in \cite{ElenaEtAlSauer2018S}. Clearly,
the associated multiresolution is significantly more economic than the
shearlet one, since $\det \Xi_j = 3^{s-1} 2$ and not $2^{2s-1}$ as in
the shearleat case.

\section{Examples}

\begin{figure}
\centering
\subfloat[$\phi_0$]{\includegraphics[width=0.3\textwidth]{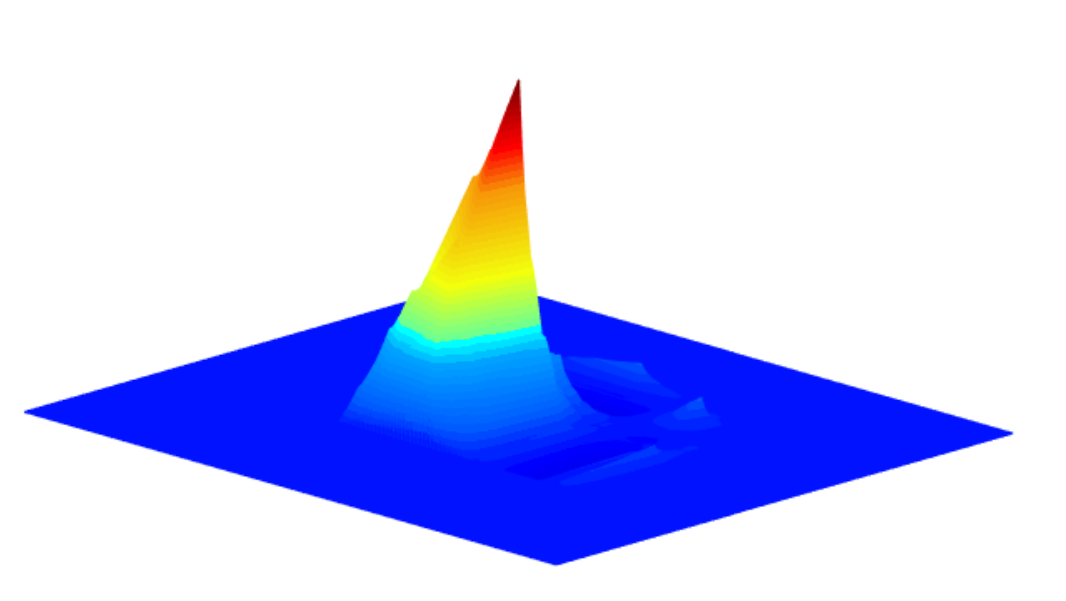}}
\subfloat[$\psi_1^0$]{\includegraphics[width=0.3\textwidth]{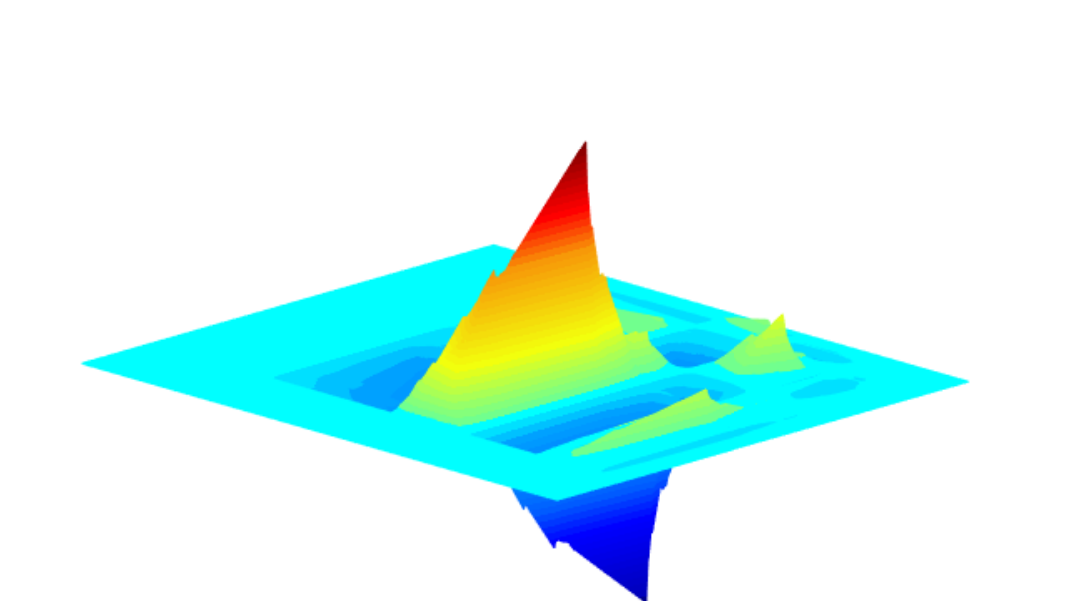}} \\
\subfloat[$\psi_2^0$]{\includegraphics[width=0.3\textwidth]{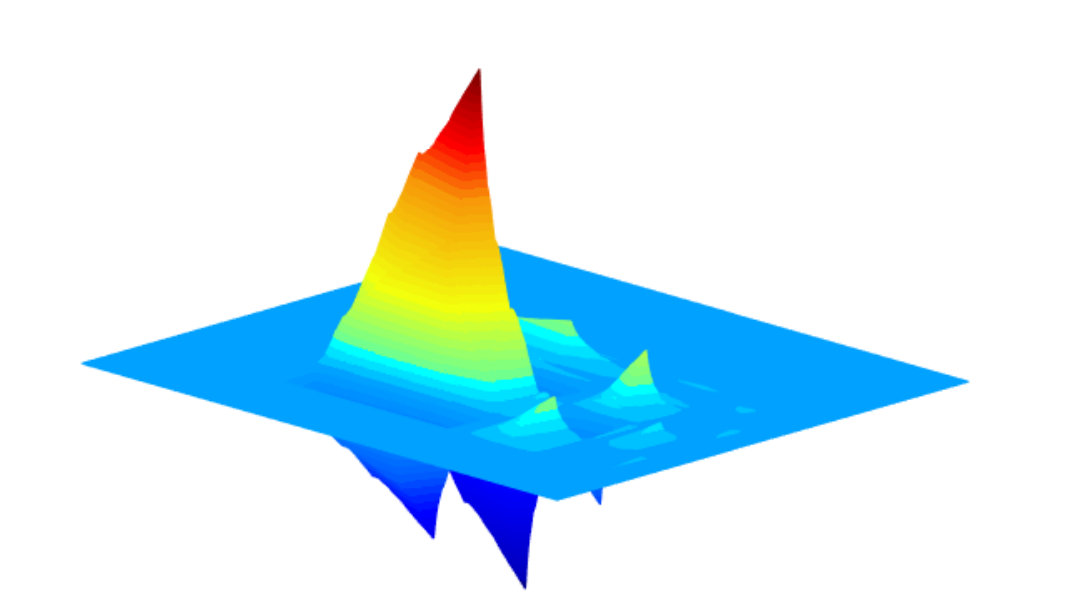}}
\subfloat[$\psi_3^0$]{\includegraphics[width=0.3\textwidth]{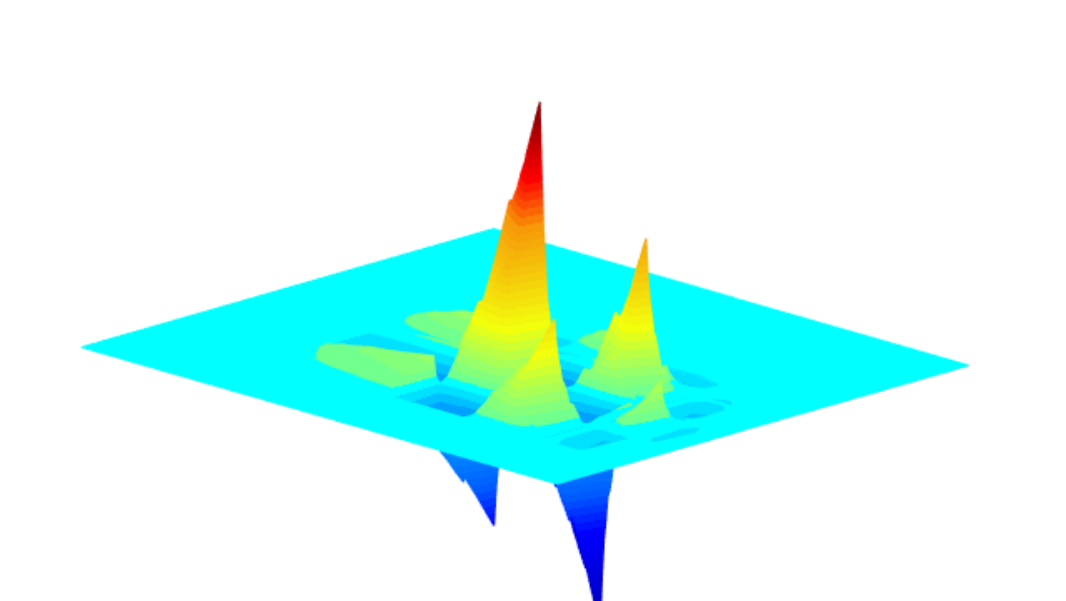}} \\
\subfloat[$\psi_4^0$]{\includegraphics[width=0.3\textwidth]{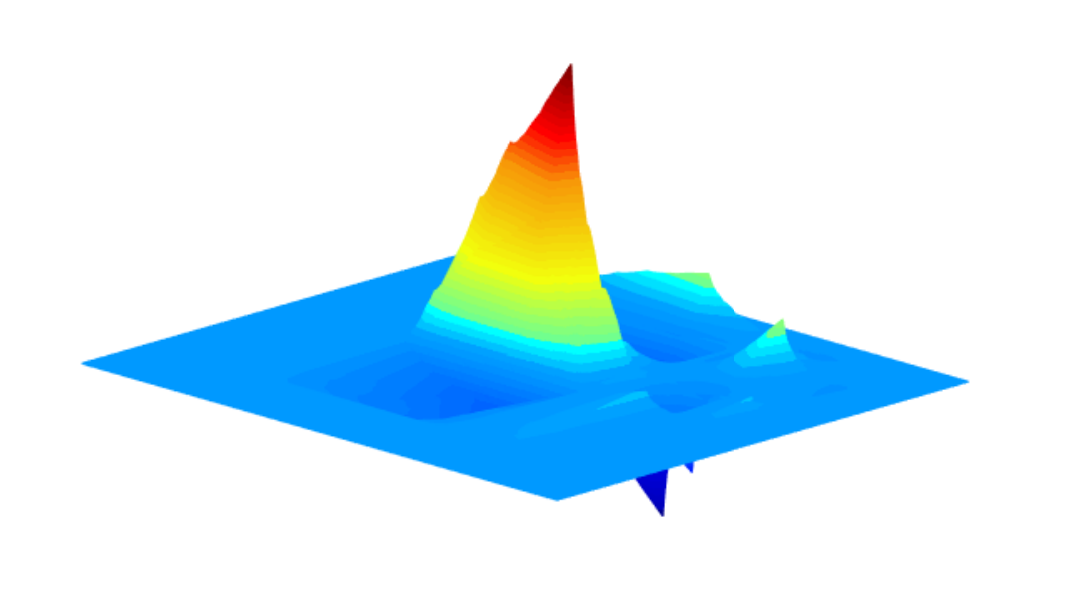}}
\subfloat[$\psi_5^0$]{\includegraphics[width=0.3\textwidth]{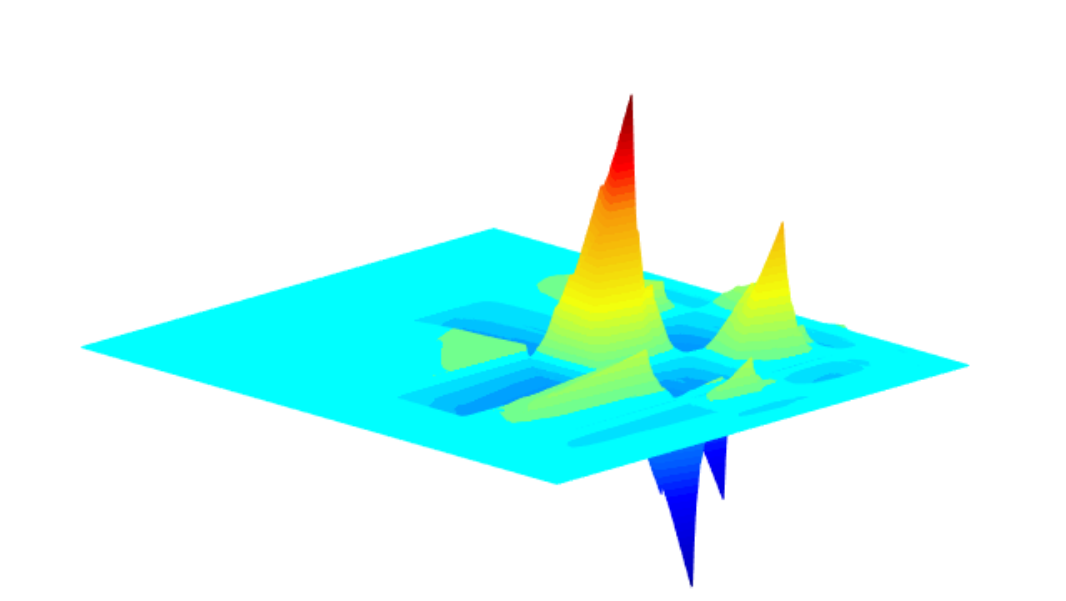}} \\
\caption{Limit functions related to the scaling matrix $\Xi_0$ in \eqref{eq:Xij_ex}.}
\label{F:limfcXi0}
\end{figure}

\begin{figure}
\centering
\subfloat[$\phi_1$]{\includegraphics[width=0.3\textwidth]{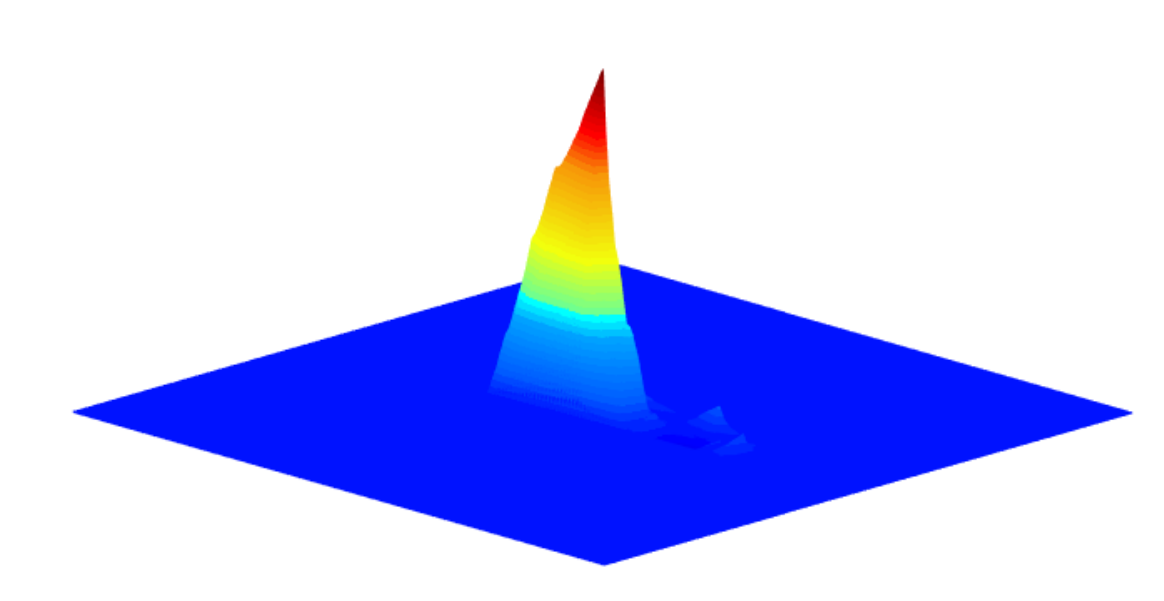}}
\subfloat[$\psi_1^1$]{\includegraphics[width=0.3\textwidth]{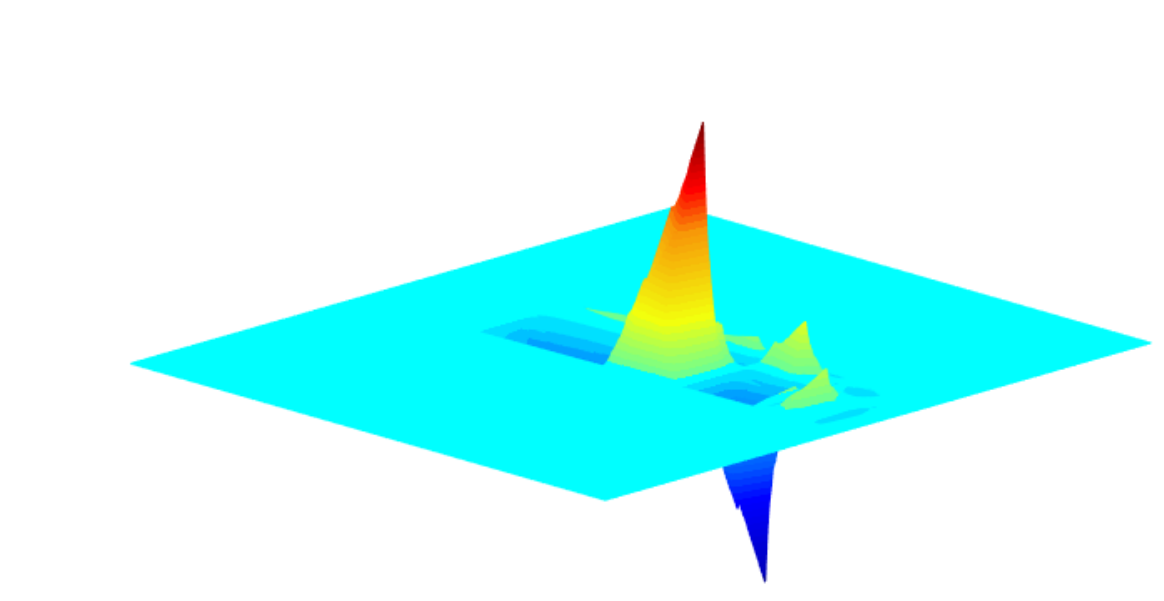}} \\
\subfloat[$\psi_2^1$]{\includegraphics[width=0.3\textwidth]{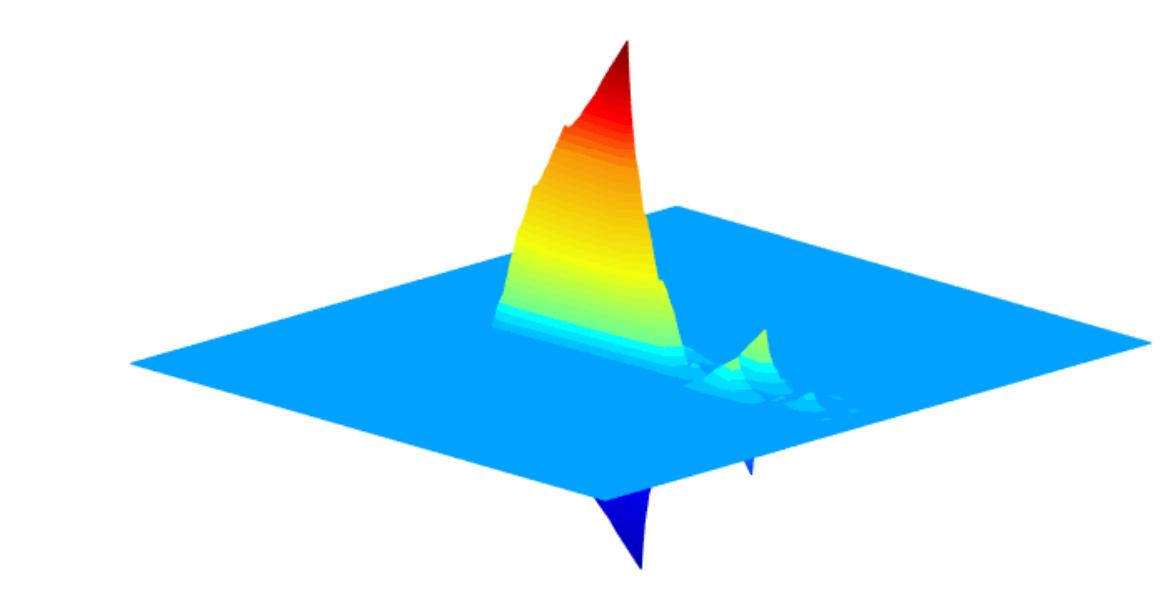}}
\subfloat[$\psi_3^1$]{\includegraphics[width=0.3\textwidth]{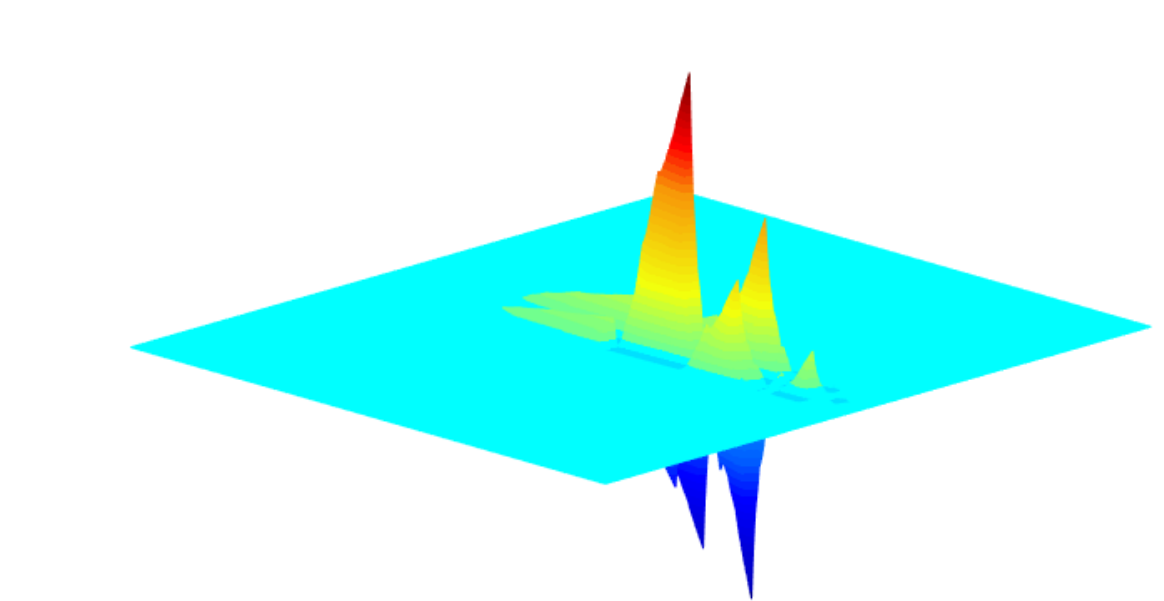}} \\
\subfloat[$\psi_4^1$]{\includegraphics[width=0.3\textwidth]{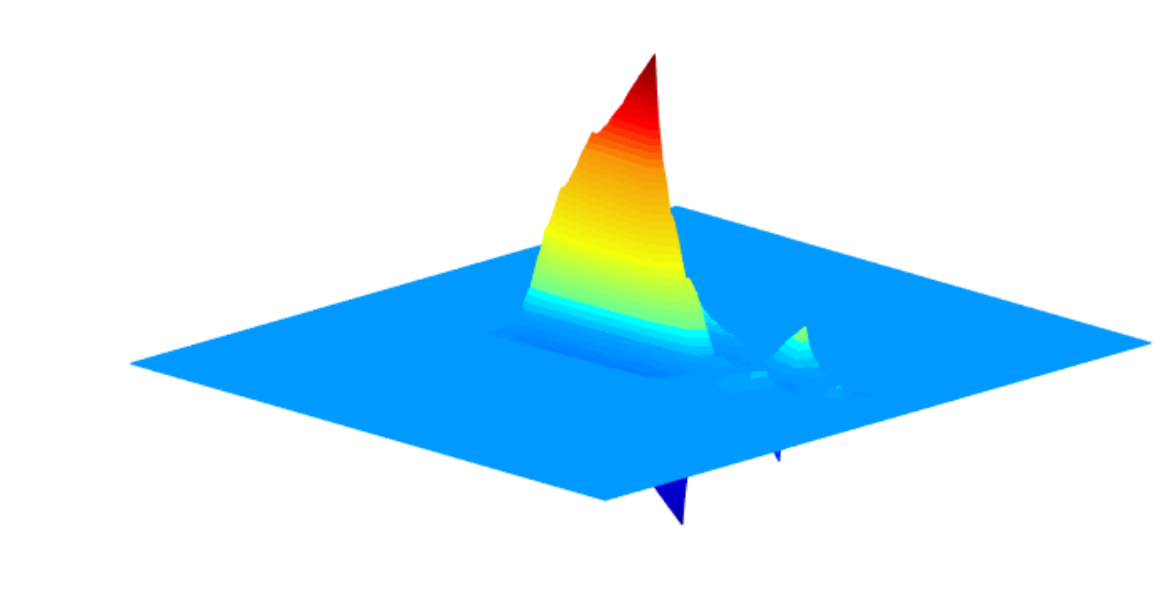}}
\subfloat[$\psi_5^1$]{\includegraphics[width=0.3\textwidth]{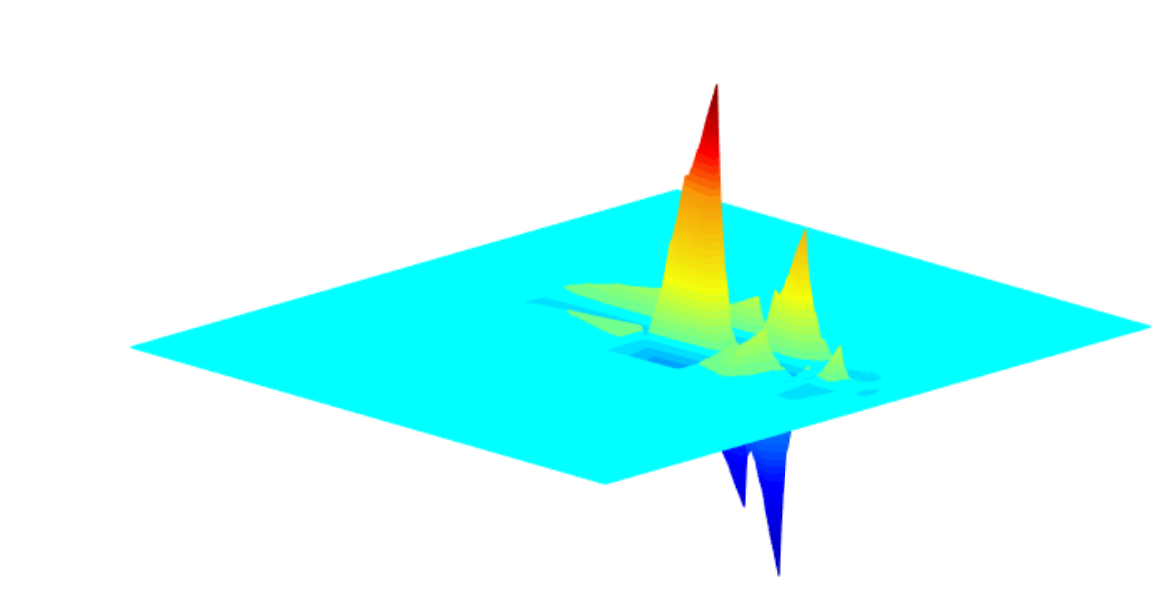}} \\
\caption{Limit functions related to the scaling matrix $\Xi_1$ in \eqref{eq:Xij_ex}.}
\label{F:limfcXi1}
\end{figure}

\noindent
Here we consider a bivariate example of matrices of the form \eqref{eq:Xij}, 
\begin{equation}
 \Xi_0=\begin{pmatrix}3 & 0\\0 & 2\end{pmatrix}, 
 \qquad \Xi_1=\begin{pmatrix} 3 & -1\\ 0 & 2\end{pmatrix},
 \label{eq:Xij_ex}
\end{equation}
that have minimum determinant possible ($\det \Xi_0=6$) because 
$\sigma_1=3$ and $\sigma_2=2$. 
The matrix $\Xi_1$ is similar to $\Xi_0$, in fact, its Smith 
factorization is
$$
 \Xi_1 =\Gamma_1^{-1}\Xi_0 \Gamma_1= 
  \begin{pmatrix}1 & 1 \\ 0 & 1\end{pmatrix} 
  \begin{pmatrix} 3 & 0 \\ 0 & 2\end{pmatrix} 
  \begin{pmatrix} 1 & -1 \\ 0 & 1\end{pmatrix} .
$$
By \eqref{eq:XiepsInvFormula} the proposed family is jointly contractive 
and it satisfies the slope resolution property, due to Theorem 
\ref{T:SlopeApproximation}.

To define an MMRA with this family of dilation matrices, we follow 
the construction of orthogonal filters described by \eqref{eq:gtensor} 
and \eqref{eq:bOrthoDef}.
For $\sigma_1=3$, we take from \cite{chui95:_const} the univariate
ternary QMF filters 
\begin{align*}
 g_1^0 &= \left( \frac{3+ \sqrt{57}}{18}, \frac{9+ \sqrt{57}}{18}, \frac{15+ \sqrt{57}}{18}, 
              \frac{15 - \sqrt{57}}{18}, \frac{9- \sqrt{57}}{18}, \frac{3- \sqrt{57}}{18} \right) ,\\
 g_1^1 &= \left( -\frac{\sqrt{2}}{2}, \sqrt{2}, -\frac{\sqrt{2}}{2} , 0, 0, 0\right) ,\\
 g_1^2 &= \frac{\sqrt{11-\sqrt{57}}}{144} \left( -21+ \sqrt{57}, -6-2\sqrt{57}, 9-5\sqrt{57}, 
                  48+8\sqrt{57}, 6+2\sqrt{57}, -36-4\sqrt{57} \right).
\end{align*}
For $\sigma_2=2$, on the other hand, we choose the univariate QMF Daubechies filters 
of order $2$
\begin{align*}
 g_2^0 &= \left( \frac{1+\sqrt{3}}{4}, \frac{3+\sqrt{3}}{4}, 
                            \frac{3-\sqrt{3}}{4}, \frac{1-\sqrt{3}}{4} \right),\\
 g_2^1 &= \left( \frac{1-\sqrt{3}}{4}, \frac{-3+\sqrt{3}}{4}, 
                            \frac{3+\sqrt{3}}{4}, \frac{-1-\sqrt{3}}{4} \right).
\end{align*}
The QMF filter system $B_0=(b_\eta^0, \, \eta \in \ZZ_6)$ with respect to the 
diagonal matrix $\Xi_0$ is composed by the tensor products
$$
  g_1^k \otimes g_2^\ell, \qquad k=0,1,2 \quad \ell=0,1.
$$
In such way, we have one low--pass filter $b_0^0$ and 5 high--pass filters 
$b_i^0,\; i=1,\ldots,5$,
\begin{equation}
 \begin{array}{cc}
   b_0^0= g_1^0 \otimes g_2^0, & b_3^0=g_1^1 \otimes g_2^1,\\
   b_1^0=g_1^0 \otimes g_2^1, & b_4^0=g_1^2 \otimes g_2^0,\\
   b_2^0=g_1^1 \otimes g_2^0, & b_5^0=g_1^2 \otimes g_2^1.
 \end{array}
 \label{eq:filtersXi0}
\end{equation}
From this family of filters we deduce the filters $B_1=(b_\eta^1, \, \eta \in \ZZ_6)$ 
associated to $\Xi_1$,
$$
  b_i^1(\cdot)=b_i^0 (\Gamma_1 \cdot), \qquad i=0,\ldots,5.
$$
By Theorem \ref{T:PsiAreWavelets} it is possible to define an orthogonal
 MMRA with dilation matrices $\Xi_0$ and $\Xi_1$. We denote with
$$
  \phi_j := \lim_{r \to \infty} S^r_{\Xi_j, b_0^j} \, \delta, \qquad j \in \ZZ_2,
$$
the scaling functions and with
$$
  \psi_i^j:= \lim_{r \to \infty} S^r_{\Xi_j, b_0^j} \, b_i^j, \qquad
  j \in \ZZ_2, \; i \in \ZZ^+_6,
$$
the wavelets functions for the prescribed family of matrices \eqref{eq:Xij_ex}.
Figures \ref{F:limfcXi0} and \ref{F:limfcXi1} depict the scaling $\phi_j$, $j \in \ZZ_2$,
and wavelet functions, $\psi_k^j$, $k \in \ZZ_6^+$, with dilation 
matrices $\Xi_0$ and $\Xi_1$, respectively.

\section*{Acknowledgement}
This research has been accomplished within Rete ITaliana di
Approssimazione (RITA). The authors Cotronei, Rossini and Volont\`e
are members of the INdAM research group GNCS.

\section*{References}


\end{document}